\documentclass{article}

\usepackage{epsfig}
\usepackage{xspace}
\usepackage{amsfonts}
\usepackage{amsmath}
\usepackage{amssymb}
\usepackage{algorithmic}
\usepackage{algorithm}
\usepackage{color}
\usepackage{harvard}

\newcommand{\oh}[1]
    {\mbox{$ {\mathcal O}( #1 ) $}}
\newcommand{\ie}
    {{\em i.e.}~}
\newcommand{\eg}
    {{\em e.g.}~}
\newcommand{\eqn}[1]
    {Equation~(\ref{eqn:#1})}

\newcommand{\alg}[1]
    {Algorithm~\ref{alg:#1}\xspace}
\begin{document}

\title{Efficient Construction,
Update and Downdate Of The Coefficients Of Interpolants Based On Polynomials Satisfying
A Three-Term Recurrence Relation}
\author{{\sc Pedro Gonnet\thanks{Corresponding author. Email: gonnet@maths.ox.ac.uk}}\\[2pt]
Departement of Computer Science, ETH Z\"urich, Switzerland and\\
Mathematical Institute, University of Oxford, United Kingdom}

\maketitle

\begin{abstract}
{In this paper, we consider methods to compute the coefficients
of interpolants relative to a basis of 
polynomials satisfying a three-term recurrence relation.
Two new algorithms are presented: the first constructs the coefficients
of the interpolation
incrementally and can be used to update the coefficients whenever
a nodes is added to or removed from the interpolation.
The second algorithm, which constructs the
interpolation coefficients by decomposing the
Vandermonde-like matrix iteratively, can not be used to update or downdate
an interpolation, yet is
more numerically stable than the first algorithm and is
more efficient when the coefficients of
multiple interpolations are to be computed
over the same set of nodes.}


\end{abstract}


\section{Introduction}

In many applications, we are interested in computing the coefficients of
a polynomial interpolation of discrete data relative to a basis
of polynomials $p_k(x)$ of increasing degree $k=0\dots n$.
In practical terms, this means that given $n+1$ function values $f_i$ at the
$n+1$ nodes $x_i$, $i=0\dots n$, we want to compute the coefficients
$c_k$, $k=0\dots n$ of a polynomial $g_n(x)$
of degree $n$ such that
\begin{equation}
    \label{eqn:problem}
    g_n(x_i) = \sum_{k=0}^n c_kp_k(x) = f_i, \quad i = 0\dots n,
\end{equation}
That is, the polynomial $g_n(x)$ interpolates the $n+1$ function values $f_i$
at the nodes $x_i$.

If we are only interested in evaluating $g_n(x)$ at different $x$,
then the method of choice is Barycentric Lagrange Interpolation
\cite{ref:Berrut2004}, which avoids representing the interpolant
in any specific base.
Such coefficient-based representations are usefull, however, if we
are interested in computing other quantities such as the
integral or derivative of the interpolant, its $L_2$-norm and/or
performing other operations on it such as transforming it 
to another interval.
Aditionally, we may also be interested in updating the coefficients
when new data is added or existing data is removed.
In \citeasnoun{ref:Gonnet2010}, such a representation is used
in an adaptive quadrature routine for just these purposes.

In the following, we will assume that the
polynomials $p_i(x)$ of degree $i$ can be constructed
using a three-term recurrence relation, which we will write as\footnote{
\citeasnoun{ref:Gautschi2004} and \citeasnoun{ref:Higham1988} use different
recurrences which can both be transformed to the representation used herein.}
\begin{equation}
    \alpha_k p_{k+1}(x) = x p_{k}(x) + \beta_k p_{k}(x) - \gamma_k p_{k-1}(x), \label{eqn:rec}
\end{equation}
with
\begin{equation*}
    p_0(x) = 1, \quad p_{-1}(x) = 0.
\end{equation*}

Examples of such polynomials are the Legendre polynomials $P_k(x)$ with
\begin{equation}
    \label{eqn:coeffs_legendre}
    \alpha_k = \frac{k}{2k-1}, \quad \beta_k = 0, \quad \gamma_k = \frac{k-1}{2k-1}
\end{equation}
or the Chebyshev polynomials $T_k(x)$ with
\begin{equation}
    \label{eqn:coeffs_chebyshev}
    \alpha_0 = 1, \quad \alpha_k = \frac{1}{2}, \quad \beta_k = 0, \quad \gamma_k = \frac{1}{2}.
\end{equation}

The coefficients
$c_i$ of \eqn{problem} can be computed solving the system of
linear equations
\begin{equation}
    \label{eqn:vandermonde-like}
    \left( \begin{array}{cccc}
        p_0(x_0) & p_1(x_0) & \dots & p_n(x_0) \\
        p_0(x_1) & p_1(x_1) & \dots & p_n(x_1) \\
        \vdots & \vdots & \ddots & \vdots \\
        p_0(x_n) & p_1(x_n) & \dots & p_n(x_n)
    \end{array} \right)
    \left( \begin{array}{c} c_0 \\ c_1 \\ \vdots \\ c_n \end{array} \right)
    =
    \left( \begin{array}{c} f_0 \\ f_1 \\ \vdots \\ f_n \end{array} \right)
\end{equation}
which can be written as
\begin{equation}
    \label{eqn:vandermonde-like2}
    \mathbf P^{(n)} \mathbf c^{(n)} = \mathbf f^{(n)}.
\end{equation}

The matrix $\mathbf P^{(n)}$ is a {\em Vandermonde-like} matrix
and the system of equations can be solved in \oh{n^3} using
Gaussian elimination.
As with the computation of the monomial coefficients, the
matrix may be ill-conditioned \cite{ref:Gautschi1983}.

\citeasnoun{ref:Bjorck1970} present an algorithm to
compute the monomial coefficients of an interpolation
without the expensive and 
potentially unstable solution of a Vandermonde
system  using Gaussian elimination, by computing first the coefficients of a Newton
interpolation and then converting these to monomial coefficients.
This approach was later extended by \citeasnoun{ref:Higham1988}
to compute the coefficients relative to any polynomial basis
satisfying a three-term recurrence relation.
Both algorithms compute the coefficients in \oh{n^2} operations
and in \citeasnoun{ref:Higham1990}, both methods are shown to
be numerically stable when a propper ordering of the nodes
$x_i$, $i=0\dots n$ is used.

In Section~\ref{sec:first}, we will re-formulate the algorithms
of Bj\"orck and Pereyra and of Higham and extend them
to update the coefficients after a {\em downdate}, \ie the removal of a node, of
an interpolation.
In Section~\ref{sec:second} we present a new algorithm for
the construction of interpolations of the type of \eqn{problem}
based on a successive decomposition of the Vandermonde-like
matrix in \eqn{vandermonde-like}.
Finally, in Section~\ref{sec:results}, we will present some
results regarding the efficiency and stability of both
algorithms.

\section{A Modification of Bj\"orck and Pereyra's and of 
Higham's Algorithms Allowing Downdates}
\label{sec:first}

\citeasnoun{ref:Bjorck1970} present an algorithm
which exploits the recursive definition of the Newton polynomials
\begin{equation}
    \label{eqn:newton_rec}
    \pi_k(x) = (x - x_{k-1}) \pi_{k-1}(x).
\end{equation}

They note that given the Newton interpolation coefficients $a_i$,
the interpolation polynomial can be constructed using Horner's scheme:
\begin{equation}
    \label{eqn:horner_q}
    q_n(x) = a_n, \quad q_k(x) = (x - x_k) q_{k+1}(x) + a_k, \quad k=n-1 \dots 0
\end{equation}
where the interpolation polynomial is $g_n(x) = q_0(x)$.

They also note that given a monomial representation for $q_k(x)$,
such as
\begin{equation*}
    q_{k}(x) = \sum_{i=0}^{n-k} b^{(k)}_i x^i
\end{equation*}
then the polynomial
$q_{k-1}(x)$ can be constructed, following the recursion in \eqn{horner_q}, as
\begin{eqnarray}
    \label{eqn:horner}
    q_{k-1}(x) & = & (x - x_{k-1}) q_{k}(x) + a_{k-1} \nonumber \\
    & = & (x - x_{k-1}) \sum_{i=0}^{n-k} b^{(k)}_i x^i + a_{k-1} \nonumber \\
    & = & \sum_{i=1}^{n-k+1} b^{(k)}_{i-1} x^i - x_{k-1} \sum_{i=0}^{n-k} b^{(k)}_i x^i + a_{k-1}  \nonumber \\
    & = & b^{(k)}_{n-k} x^{n-k+1} + \sum_{i=1}^{n-k} (b^{(k)}_{i-1} - x_{k-1}b^{(k)}_i) x^i + -b^{(k)}_0x_{k-1} + a_{k-1}.
\end{eqnarray}
From \eqn{horner} we can then extract the new coefficients $b^{(k-1)}_i$:
\begin{equation}
    \label{eqn:b_update}
    b^{(k-1)}_i = \left\{ \begin{array}{ll}
        b^{(k)}_{i-1}, & i = n-k+1, \\
        b^{(k)}_{i-1} - x_{k-1} b^{(k)}_i,\quad & 1 \leq i \leq n-k, \\
        -b^{(k)}_0x_{k-1} + a_{k-1}, & i = 0. \\
    \end{array} \right.
\end{equation}

\citeasnoun{ref:Higham1988} uses the same approach, yet represents
the Newton polynomials as a linear combination of polynomials satisfying
a three-term recurrence relation.
Using such a representation
\begin{equation}
    \label{eqn:higham_qk}
    q_k(x) = \sum_{i=0}^{n-k} c^{(k)}_i p_i(x)
\end{equation}
he computes $q_{k-1}(x)$ by expanding the recursion in \eqn{horner_q}
using the representation \eqn{higham_qk}:
\begin{eqnarray}
    \label{eqn:horner_expand1}
    q_{k-1}(x) & = & (x - x_{k-1}) \sum_{i=0}^{n-k} c^{(k)}_i p_i(x) + a_{k-1} \nonumber \\
    & = & \sum_{i=0}^{n-k} c^{(k)}_i xp_i(x) - x_{k-1} \sum_{i=0}^{n-k} c^{(k)}_i p_i(x) + a_{k-1} \nonumber \\
    & = & \sum_{i=0}^{n-k} c^{(k)}_i \left( \alpha_i p_{i+1}(x) - \beta_i p_i(x) + \gamma_i p_{i-1}(x) \right) - x_{k-1} \sum_{i=0}^{n-k} c^{(k)}_i p_i(x) + a_{k-1}.
\end{eqnarray}
Expanding \eqn{horner_expand1} for the individual $p_k(x)$, and
keeping in mind that $p_{-1}(x) = 0$, we obtain
\begin{eqnarray}
    \label{eqn:horner_expand2}
    q_{k-1}(x) & = & \sum_{i=1}^{n-k+1} c^{(k)}_{i-1} \alpha_{i-1} p_{i}(x) - \sum_{i=0}^{n-k} c^{(k)}_i \left( x_{k-1} + \beta_i \right) p_i(x) \nonumber \\
        & & + \sum_{i=0}^{n-k-1} c^{(k)}_{i+1} \gamma_{i+1} p_{i}(x) + a_{k-1}.
\end{eqnarray}
By shifting the sums in \eqn{horner_expand2} and re-grouping around
the individual $p_k(x)$ we finally obtain
\begin{eqnarray}
    \label{eqn:horner2}
    q_{k-1}(x) & = & c^{(k)}_{n-k}\alpha_{n-k}p_{n-k+1}(x) + \left( c^{(k)}_{n-k-1}\alpha_{n-k-1} - c^{(k)}_{n-k}(x_{k-1} + \beta_{n-k})\right)p_{n-k}(x) \nonumber \\
    &   & + \sum_{i=1}^{n-k-1} \left(c^{(k)}_{i-1} \alpha_{i-1} - c^{(k)}_i (x_{k-1} + \beta_i) + c^{(k)}_{i+1} \gamma_{i+1} \right)p_i(x) \nonumber \\
    &   & - c^{(k)}_0(x_{k-1}+\beta_0) + c^{(k)}_1 \gamma_1 + a_{k-1}.
\end{eqnarray}
Higham then extracts the new coefficients $c^{(k-1)}_i$ from
\eqn{horner2} as:
\begin{equation}
    \label{eqn:c_update}
    c^{(k-1)}_i = \left\{ \begin{array}{ll}
        c^{(k)}_{i-1} \alpha_{i-1}, \quad & i = n-k+1, \\
        c^{(k)}_{i-1} \alpha_{i-1} - c^{(k)}_{i}(x_{k-1} + \beta_{i}), \quad & i = n-k, \\
        c^{(k)}_{i-1} \alpha_{i-1} - c^{(k)}_i (x_{k-1} + \beta_i) + c^{(k)}_{i+1} \gamma_{i+1}, \quad & 1 \leq i < n-k, \\
        -c^{(k)}_0(x_{k-1}+\beta_0) + c^{(k)}_1 \gamma_1 + a_{k-1}, & i = 0
    \end{array} \right.
\end{equation}
In both algorithms, the interpolating polynomial is constructed by
first computing the divided differences
\begin{equation}
    \label{eqn:newton_coeffs}
    a_i = f[x_0,\dots,x_i],\quad i=0\dots n
\end{equation}
and, starting with
$q_n(x)=a_n$, and hence $c^{(n)}_0=a_n$ or $b^{(n)}_0=a_n$,
successively updating the coefficients per \eqn{c_update} or \eqn{b_update}
respectively.

Alternatively, we could use the same approach to
compute the coefficients of the Newton polynomials themselves
\begin{equation*}
    \pi_k(x) = \sum_{i=0}^{k} \eta^{(k)}_i p_i(x).
\end{equation*}

Expanding the recurrence relation in \eqn{newton_rec} analogously to
\eqn{horner2}, we get
\begin{eqnarray}
    \pi_{k+1}(x) & = & \eta^{(k)}_{k}\alpha_{k}p_{k+1}(x) + \left( \eta^{(k)}_{k-1}\alpha_{k-1} - \eta^{(k)}_{k}(x_k + \beta_{k}) \right)p_k(x) \nonumber \\
    & & + \sum_{i=1}^{k-1} \left( \eta^{(k)}_{i-1} \alpha_{i-1} - \eta^{(k)}_i (x_k + \beta_i) + \eta^{(k)}_{i+1} \gamma_{i+1} \right) p_i(x) \nonumber \\
    & & - \eta^{(k)}_0(x_k+\beta_0) + \eta^{(k)}_1 \gamma_1.
\end{eqnarray}

We initialize with $\eta^{(0)}_0 = 1$ and use
\begin{equation}
    \label{eqn:eta_rec}
    \eta^{(k+1)}_i = \left\{ \begin{array}{ll}
        \eta^{(k)}_{i-1}\alpha_{i-1}, \quad & i = k+1, \\
        \eta^{(k)}_{i-1}\alpha_{i-1} - \eta^{(k)}_{i}(x_k + \beta_{i}), \quad & i = k, \\
        \eta^{(k)}_{i-1} \alpha_{i-1} - \eta^{(k)}_i (x_k + \beta_i) + \eta^{(k)}_{i+1} \gamma_{i+1}, \quad & 1 \leq i < k, \\
        -\eta^{(k)}_0(x_k+\beta_0) + \eta^{(k)}_1 \gamma_1, & i = 0,
    \end{array} \right.
\end{equation}
to compute the coefficients for $\pi_k(x)$, $k=1\dots n$.
Alongside this computation, we can also compute the coefficients of
a sequence of polynomials $g_k(x)$ of increasing degree $k$
\begin{equation*}
    g_k(x) = \sum_{i=0}^k c^{(k)}_i p_i(x)
\end{equation*}
initializing with $c^{(0)}_0 = a_0$, where the $a_i$ are still the 
Newton coefficients computed and used above.
The subsequent coefficients $c^{(k)}_i$, $k=1\dots n$ are computed
using
\begin{equation}
    \label{eqn:c_rec}
    c^{(k)}_i = \left\{ \begin{array}{ll}
        \eta^{(k)}_i a_k, \quad & i = k, \\
        c^{(k-1)}_i + \eta^{(k)}_i a_k, \quad & 0 \leq i < k.
    \end{array} \right.
\end{equation}

This {\em incremental} construction of the coefficients, which is 
equivalent to effecting the summation of the weighted Newton polynomials
and is referred
to by Bj\"orck and Pereyra as the ``progressive algorithm'', can be
used to efficiently update an interpolation.
If the coefficients $\eta^{(n)}_i$ and $c^{(n)}_i$ are stored and
a new node $x_{n+1}$ and function value $f_{n+1}$ are added
to the data, a new coefficient $a_{n+1}$ can be computed per
\eqn{newton_coeffs}, the coefficients $\eta^{(n+1)}_i$ computed
per \eqn{eta_rec} and, finally, the $c^{(n)}_i$ updated per
\eqn{c_rec}, resulting in the coefficients $c^{(n+1)}_i$ for the updated 
interpolation polynomial $g_{n+1}(x)$.

\begin{algorithm}
    \caption{Incremental construction of $g_n(x)$}
    \label{alg:incr}
    \begin{algorithmic}[1]
        \STATE $c^{(0)}_0 \leftarrow f_0$ \hfill ({\em init $\mathbf c^{(0)}$})
        \STATE $\eta^{(1)}_0 \leftarrow -x_0 - \beta_0$, $\eta^{(1)}_1 \leftarrow \alpha_0$ \hfill ({\em init $\boldsymbol{\eta}^{(1)}$})
        \FOR{$k=1 \dots n$}
            \STATE $v_0 \leftarrow 0$, $v_1 \leftarrow x_k$ \hfill ({\em init $\mathbf v$})
            \FOR{$i = 2 \dots k$}
                \STATE $v_i \leftarrow \left( ( x_k + \beta_{i-1}) v_{i-1} - \gamma_{i-1}v_{i-2} \right) / \alpha_{i-1}$ \hfill ({\em compute the $v_i$})
            \ENDFOR
            \STATE $g_k \leftarrow \mathbf v(0:k-1)^\mathsf{T} \mathbf c^{(k-1)}$ \hfill ({\em compute $g_{k-1}(x_k)$})
            \STATE $\pi_k \leftarrow \mathbf v^\mathsf{T} \mathbf \eta^{(k)}$ \hfill ({\em compute $\pi_k(x_k)$})
            \STATE $a_k \leftarrow ( f_k - g_k ) / \pi_k$ \hfill ({\em compute $a_k$, \eqn{anp1}})
            \STATE $\mathbf c^{(k)} \leftarrow [ \mathbf c^{(k-1)} ; 0 ] + a_k \boldsymbol{\eta}^{(k)}$ \hfill ({\em compute the new $\mathbf c^{(k)}$, \eqn{c_rec}})
            \STATE $\boldsymbol{\eta}^{(k+1)} \leftarrow [ 0 ; \underline{\alpha}(0:k) .* \boldsymbol{\eta}^{(k)} ] - [ (x_k + \underline{\beta}(0:k)) .* \boldsymbol{\eta}^{(k)} ; 0 ] +$ \\ \hspace{\algorithmicindent} $[ \underline{\gamma}(1:k) .* \boldsymbol{\eta}^{(k)}(1:k) ; 0 ; 0 ]$ \hfill ({\em compute the new $\boldsymbol{\eta}^{(k+1)}$, \eqn{eta_rec}})
        \ENDFOR
    \end{algorithmic}
\end{algorithm}

We can re-write the recursion for the coefficients $\eta^{(k)}_i$ 
of the Newton polynomials in matrix-vector notation as
\begin{equation*}
    \left( \mathbf T^{(k+1)} - \underline{\mathbf I}_0 x_{k} \right) \boldsymbol{\eta}^{(k)} = \boldsymbol{\eta}^{(k+1)}
\end{equation*}
where $\mathbf T^{(k+1)}$ is the $(k+2)\times (k+1)$ tri-diagonal matrix
\begin{equation*}
    \mathbf T^{(k+1)} = \left( \begin{array}{ccccc}
        -\beta_0 & \gamma_1 \\
        \alpha_0 & -\beta_1 & \gamma_2 \\
        & \ddots & \ddots & \ddots \\
        & & \alpha_{k-2} & -\beta_{k-1} & \gamma_k \\
        & & & \alpha_{k-1} & -\beta_k \\
        & & & & \alpha_k
    \end{array} \right)
\end{equation*}
and $\underline{\mathbf I}_0 x_{k}$ is a $(k+2)\times (k+1)$ 
matrix with $x_{k}$ in the diagonal and zeros elsewhere
\begin{equation*}
    \underline{\mathbf I}_0 x_{k} = \left( \begin{array}{cccc}
        x_{k} \\
        & x_{k} \\
        & & \ddots \\
        & & & x_{k} \\
        0 & 0 & \dots & 0
    \end{array} \right).
\end{equation*}
The vectors $\boldsymbol{\eta}^{(k)} = ( \eta^{(k)}_0 , \eta^{(k)}_1 , \dots , \eta^{(k)}_{k} )^\mathsf{T}$
and $\boldsymbol{\eta}^{(k+1)} = ( \eta^{(k+1)}_0 , \eta^{(k+1)}_1 , \dots , \eta^{(k+1)}_{k+1} )^\mathsf{T}$
contain the coefficients of the $k$th and $(k+1)$st Newton polynomial
respectively.

Given the vector of coefficients $\mathbf c^{(n)} = ( c^{(n)}_0 , c^{(n)}_1 , \dots , c^{(n)}_n )^\mathsf{T}$
of an interpolation polynomial $g_n(x)$ of 
degree $n$ and the vector of coefficients $\boldsymbol{\eta}^{(n+1)}$ of
the $(n+1)$st Newton polynomial over the $n+1$ nodes, we can update the
interpolation for a new node $x_{n+1}$ and function value $f_{n+1}$
as follows:
Instead of computing the new Newton interpolation coefficient $a_{n+1}$ using
the divided
differences as in \eqn{newton_coeffs}, we choose $a_{n+1}$ such that
the new interpolation constraint
\begin{equation*}
    g_{n+1}(x_{n+1}) = g_n(x_{n+1}) + a_{n+1} \pi_{n+1}(x_{n+1}) = f_{n+1}
\end{equation*}
is satisfied, resulting in 
\begin{equation}
    \label{eqn:anp1}
    a_{n+1} = \frac{f_{n+1} - g_n(x_{n+1})}{\pi_{n+1}(x_{n+1})}
\end{equation}
which can be computed by evaluating $g_n(x_{n+1})$ and $\pi_{n+1}(x_{n+1})$.
Note that since $\pi_{n+1}(x_i)=0$ for $i=0\dots n$, the addition of any
multiple of $\pi_{n+1}(x)$ to $g_n(x)$ does not affect the interpolation
at the other nodes at all.
This expression for $a_{n+1}$ is used instead of the divided difference
since we have not explicitly stored
the previous $a_i$, $i = 0\dots n$, which are needed for the recursive
computation of the latter.

We then update the coefficients of the interpolating polynomial
using
\begin{equation*}
    \mathbf c^{(n+1)} = \left( \begin{array}{c} \mathbf c^{(n)} \\ 0 \end{array} \right)
        + a_{n+1} \boldsymbol{\eta}^{(n+1)}
\end{equation*}
and then the coefficients of the Newton polynomial using
\begin{equation*}
     \boldsymbol{\eta}^{(n+2)} = \left( \mathbf T^{(n+2)} - \underline{\mathbf I}_0 x_{n+1} \right) \boldsymbol{\eta}^{(n+1)}
\end{equation*}
such that it is ready for further updates.
Starting with $\eta^{(0)}_0 = 1$ and $n=0$, this update can be used
to construct $g_n(x)$ by adding each $x_i$ and $f_i$,
$i=0 \dots n$, successively.

The complete algorithm doing just that is shown in 
Algorithm~\ref{alg:incr}.
The addition of each $n$th node requires \oh{n} operations,
resulting in a total of \oh{n^2} operations for the construction
of an $n$-node interpolation.

This is essentially the progressive algorithm of Bj\"orck and
Pereyra, yet instead of storing the Newton coefficients $a_i$,
we store the coefficients $\eta^{(n+1)}_i$ of the last Newton polynomial.
This new representation offers no obvious advantage for the 
update, other than that it can be easily {\em reversed}:

Given an interpolation over a set of $n+1$ nodes $x_i$ and 
function values $f_i$, $i=0\dots n$ defined by the coefficients
$c^{(n)}_i$ and given the coefficients $\eta^{(n+1)}_i$ of the $(n+1)st$ Newton
polynomial over the same nodes, we will {\em downdate} the
interpolation by removing the function value $f_j$ at the node $x_j$.
The resulting polynomial of degree $n-1$ will still interpolate the
remaining $n$ nodes.

We start by removing the root $x_j$ from the $(n+1)$st Newton
polynomial by solving
\begin{equation*}
     \left( \mathbf T^{(n+1)} - \underline{\mathbf I}_0 x_j \right) \boldsymbol{\eta}^{(n)} = \boldsymbol{\eta}^{(n+1)}
\end{equation*}
for the vector of coefficients $\boldsymbol{\eta}^{(n)}$.
Since $x_j$ is a root of $\pi_{n+1}(x)$, the system is over-determined 
yet has a unique solution\footnote{Note that the $n \times (n+1)$ matrix 
$\left(\mathbf T^{(n+1)}-\underline{\mathbf I}_0x_j\right)^\mathsf{T}$ has rank $n$
and the null space $\mathbf p(x_j) = \left( p_0(x_j), p_1(x_j) , \dots , p_{n+1}(x_j) \right)^\mathsf{T}$
since for $\mathbf v = \left(\mathbf T^{(n+1)}-\underline{\mathbf I}_0x_j\right)^\mathsf{T} \mathbf p(x_j)$,
$v_i = \alpha_i p_{i+1}(x_j) - (x_j + \beta_i)p_i(x_j) + \gamma_i p_{i-1}(x_j) = 0$
by the definition in \eqn{rec} and the right-hand side $\boldsymbol{\eta}^{(n+1)}$
is consistent.}.
We can therefore remove the first row of $(\mathbf T^{(n+1)}- \underline{\mathbf I}_0 x_j)$
and the first entry of $\boldsymbol{\eta}^{(n+1)}$, resulting in
the upper-tridiagonal system of linear equations
\begin{equation}
    \label{eqn:second_linsys}
    \left( \begin{array}{ccccc}
        \alpha_0 & -(x_j+\beta_1) & \gamma_2 \\
        & \ddots & \ddots & \ddots \\
        & & \alpha_{n-2} & -(x_j+\beta_{n-1}) & \gamma_n \\
        & & & \alpha_{n-1} & -(x_j+\beta_n) \\
        & & & & \alpha_n
    \end{array} \right)
    \left( \begin{array}{c} \eta^{(n)}_0 \\ \eta^{(n)}_1 \\ \vdots \\ \eta^{(n)}_n \end{array} \right)
    =
    \left( \begin{array}{c} \eta^{(n+1)}_1 \\ \eta^{(n+1)}_2 \\ \vdots \\ \eta^{(n+1)}_{n+1} \end{array} \right)
\end{equation}
which can be conveniently solved in \oh{n} using back-substitution.

Once we have our downdated $\boldsymbol{\eta}^{(n)}$, and thus
the downdated Newton polynomial $\pi_n(x)$, we can
downdate the coefficients of $g_n(x)$ by computing
\begin{equation*}
    g_{n-1}(x) = g_n(x) - a^{\star}_j \pi_n(x)
\end{equation*}
where the Newton coefficient $a^{\star}_j$ would need to
be re-computed from the divided difference over all nodes 
{\em except} $x_j$.
We can avoid this computation by noting that $g_{n-1}(x)$ has to be
of degree $n-1$ and therefore the highest coefficient of $g_n(x)$,
$c^{(n)}_n$, must disappear.
This is the case when
\begin{equation*}
    c^{(n-1)}_n = c^{(n)}_n - a^{\star}_j \eta^{(n)}_n = 0
\end{equation*}
and therefore
\begin{equation*}
    a^{\star}_j = \frac{c^{(n)}_n}{\eta^{(n)}_n}.
\end{equation*}

Using this $a^{\star}_j$, we can the compute the coefficients of $g_{n-1}(x)$
as
\begin{equation}
    \label{eqn:first_down}
    c^{(n-1)}_i = c^{(n)}_i - \frac{c^{(n)}_n}{\eta^{(n)}_n} \eta^{(n)}_i, \quad i=1 \dots n-1.
\end{equation}

The whole process is shown in Algorithm~\ref{alg:downdate}.
The downdate of an $n$-node interpolation requires \oh{n}
operations.

\begin{algorithm}
    \caption{Remove a function value $f_j$ at the node $x_j$ from the interpolation
        given by the coefficients $\mathbf c^{(n)}$}
    \label{alg:downdate}
    \begin{algorithmic}[1]
        \STATE $\eta^{(n)}_n \leftarrow \eta^{(n+1)}_{n+1} / \alpha_n$ 
            \hfill ({\em compute $\boldsymbol{\eta}^{(n)}$ from $\boldsymbol{\eta}^{(n+1)}$ using back-substitution})
        \STATE $\eta^{(n)}_{n-1} \leftarrow \left( \eta^{(n+1)}_n + (x_j + \beta_n)\eta^{(n)}_n \right) / \alpha_{n-1}$
        \FOR{$i=n-2 \dots 0$}
            \STATE $\eta^{(n)}_i \leftarrow \left( \eta^{(n+1)}_{i+1} + (x_j + \beta_{i+1})\eta^{(n)}_{i+1} - \gamma_{i+2}\eta^{(n)}_{i+2} \right) / \alpha_i$
        \ENDFOR
        \STATE $a_j \leftarrow c^{(n)}_n / \eta^{(n)}_n$ 
            \hfill ({\em compute the coefficient $a_j$})
        \STATE $\mathbf c^{(n-1)} \leftarrow \mathbf c^{(n)} - a_j \boldsymbol{\eta}^{(n)}$ 
            \hfill ({\em compute the new coefficients $\mathbf c^{(n-1)}$})\label{alg:downdate_last}
    \end{algorithmic}
\end{algorithm}

\section{A New Algorithm for the Construction of Interpolations}
\label{sec:second}

Returning to the representation in \eqn{vandermonde-like}, we can
try to solve the Vandermonde-like system of linear equations directly.
The matrix has some special characteristics which we can exploit
to achieve better performance and stability than when
using Gaussian elimination or even the algorithms of Bj\"orck and Pereyra,
Higham or the one described in the previous section.

We start by de-composing the $(n+1)\times(n+1)$ Vandermonde-like
matrix $\mathbf P^{(n)}$ as follows:
\begin{equation*}
    \mathbf P^{(n)} = \left( \begin{array}{c|c}
        \quad  \rule[-0.9cm]{0cm}{2cm} \mathbf P^{(n-1)} \quad  & \mathbf p^{(n)}  \\ \hline
        \mathbf q^\mathsf{T} & p_n(x_n)
    \end{array} \right).
\end{equation*}
The sub-matrix $\mathbf P^{(n-1)}$ is a Vandermonde-like matrix
analogous to $\mathbf P^{(n)}$.
The column $\mathbf p^{(n)}$ contains the $n$th polynomial evaluated
at the nodes $x_i$, $i=0 \dots n-1$
\begin{equation*}
    \mathbf p^{(n)} = \left( \begin{array}{c} p_n(x_0) \\ p_n(x_1) \\ \vdots \\ p_n(x_{n-1}) \end{array} \right)
\end{equation*}
and the vector $\mathbf q^\mathsf{T}$ contains the values of the first $0\dots n-1$
polynomials at the node $x_n$
\begin{equation*}
    \mathbf q^\mathsf{T} = \left( p_0(x_n) , p_1(x_n) , \dots , p_{n-1}(x_n) \right).
\end{equation*}

Inserting this into the product in \eqn{vandermonde-like2}, we obtain
\begin{equation*}
    \left( \begin{array}{c|c}
        \quad  \rule[-0.9cm]{0cm}{2cm} \mathbf P^{(n-1)} \quad  & \mathbf p^{(n)}  \\ \hline
        \mathbf q^\mathsf{T} & p_n(x_n)
    \end{array} \right)
    \left( \begin{array}{c} \rule[-0.9cm]{0cm}{2cm} \mathbf c^{(n-1)} \\ \hline c_n \end{array} \right)
    =
    \left( \begin{array}{c} \rule[-0.9cm]{0cm}{2cm} \mathbf f^{(n-1)} \\ \hline f_n \end{array} \right) \nonumber \\
\end{equation*}
which, when effected, results in the pair of equations
\begin{eqnarray}
    \label{eqn:second_two}
    \mathbf P^{(n-1)} \mathbf c^{(n-1)} + \mathbf p^{(n)} c_n & = & \mathbf f^{(n-1)} \\
    \mathbf q^\mathsf{T} \mathbf c^{(n-1)} + p_n(x_n) c_n & = & f_n, \nonumber
\end{eqnarray}
where the vectors $\mathbf c^{(n-1)} = \left( c_0 , c_1 , \dots , c_{n-1} \right)^\mathsf{T}$
and $\mathbf f^{(n-1)} = \left( f_0 , f_1 , \dots , f_{n-1} \right)^\mathsf{T}$
contain the first $n$ coefficients or function values respectively.

Before trying to solve \eqn{second_two}, we note that the columns of
the matrix
$\mathbf P^{(n-1)}$ contain the first $0 \dots n-1$ polynomials
evaluated at the same $n$ nodes each.
Similarly, $\mathbf q^\mathsf{T}$ contains the same polynomials
evaluated at the node $x_n$.
Since the polynomials in the columns are of degree $< n$ and they are evaluated
at $n$ points, $\mathbf P^{(n-1)}$ actually contains enough data
to extrapolate the values of these polynomials at $x_n$.
Using Lagrange interpolation we can write
\begin{equation}
    \label{eqn:second_qi}
    q_i = \sum_{j=0}^{n-1} \ell_j^{(n)}(x_n) P^{(n-1)}_{j,i}
\end{equation}
where the
\begin{equation}
    \label{eqn:lagrange_poly}
    \ell^{(n)}_j(x) = \prod_{\substack{j=0\\j \ne i}}^(n-1) \frac{x - x_j}{x_i - x_j}
\end{equation}
are the Lagrange polynomials over the first $n$
nodes $x_i$, $i=0 \dots n-1$.
We can write \eqn{second_qi} as
\begin{equation}
    \label{eqn:second_qt}
    \mathbf q^\mathsf{T} = \boldsymbol{\ell}^{(n)} \mathbf P^{(n-1)}
\end{equation}
where the entries of the $1 \times n$ vector $\boldsymbol{\ell}^{(n)}$ are
\begin{equation*}
    \ell^{(n)}_i = \ell^{(n)}_i(x_n).
\end{equation*}

The entries of $\boldsymbol{\ell}^{(n)}$ can be computed recursively.
Using the definition in \eqn{lagrange_poly}, we define 
\begin{equation*}
    w_n = \prod_{j=0}^{n-1} (x_n - x_j).
\end{equation*}
and re-write $\ell^{(n)}_i$ as
\begin{equation}
    \label{eqn:ellni1}
    \ell^{(n)}_i = \frac{w_n}{x_n - x_i} \left[ \prod^{n-1}_{\substack{j=0\\j \ne i}} (x_i - x_j) \right]^{-1}.
\end{equation}

Using the previous $\ell^{(n-1)}_i$ and $w_{n-1}$, we can re-write \eqn{ellni1}
as
\begin{equation*}
    \ell^{(n-1)}_i = \frac{\ell^{(n-1)}_i}{w_{n-1}}(x_{n-1}-x_i)
\end{equation*}
which, re-inserted into \eqn{ellni1}, gives
\begin{equation}
    \ell^{(n)}_i =
        \frac{w_n}{x_n - x_i} \frac{\ell^{(n-1)}_i}{w_{n-1}} \frac{x_{n-1}-x_i}{x_i - x_{n-1}} =
        - \frac{\ell^{(n-1)}_i}{x_n - x_i} \frac{w_n}{w_{n-1}} \label{eqn:ellni}
\end{equation}
for all $i < n-1$.
We then compute the last entry $i = n-1$ using
\begin{equation}
    \label{eqn:ellnnm1}
    \ell^{(n)}_{n-1} = \frac{w_n}{w_{n-1}(x_{n} - x_{n-1})}.
\end{equation}
Therefore, starting with $\ell^{(1)}_0 = 1$, we can construct all
the $\boldsymbol{\ell}^{(k)}$, $k=2 \dots n$ successively.
Since, given $\boldsymbol{\ell}^{(k)}$, $k=1\dots n$ the construction 
of each additional $\boldsymbol{\ell}^{(n+1)}$
requires \oh{n} operations, the construction of all the
$\boldsymbol{\ell}^{(k)}$, $k=1\dots n$ requires a total
of \oh{n^2} operations.

Returning to the Vandermonde-like matrix,
inserting \eqn{second_qt} into \eqn{second_two}, we obtain
\begin{eqnarray}
    \label{eqn:second_two2}
    \mathbf P^{(n-1)} \mathbf c^{(n-1)} + \mathbf p^{(n)} c_n & = & \mathbf f^{(n-1)} \nonumber \\
    \boldsymbol{\ell}^{(n)} \mathbf P^{(n-1)}  \mathbf c^{(n-1)} + p_n(x_n) c_n & = & f_n.
\end{eqnarray}
Multiplying the first line in \eqn{second_two2} with $\boldsymbol{\ell}^{(n)}$
from the left and subtracting the bottom equation from the top one we obtain
\begin{equation*}
    \boldsymbol{\ell}^{(n)} \mathbf p^{(n)} c_n - p_n(x_n) c_n = \boldsymbol{\ell}^{(n)} \mathbf f^{(n-1)} - f_n
\end{equation*}
from which we can finally isolate the coefficient $c_n$:
\begin{equation}
    \label{eqn:cn}
    c_n = \frac{\boldsymbol{\ell}^{(n)} \mathbf f^{(n-1)} - f_n}{\boldsymbol{\ell}^{(n)} \mathbf p^{(n)} - p_n(x_n)}.
\end{equation}
Having computed $c_n$, we can now re-insert it into the first Equation
in \eqn{second_two2}, resulting in the new system
\begin{equation}
    \label{eqn:second_new}
    \mathbf P^{(n-1)} \mathbf c^{(n-1)} = \mathbf f^{(n-1)} - \mathbf p^{(n)}c_n
\end{equation}
in the remaining coefficients $\mathbf c^{(n-1)}$.

Applying this computation recursively to $\mathbf P^{(n)}$, $\mathbf P^{(n-1)}$,
$\dots$, $\mathbf P^{(1)}$ we can compute the interpolation coefficients
$\mathbf c^{(n)}$.
The final coefficient $c_0$ can be computed as
\begin{equation*}
    c_0 = f_0 / P^{(0)}_{1,1}.
\end{equation*}
The complete algorithm is shown in \alg{direct}.
Since the construction of the $\boldsymbol{\ell}^{(k)}$, 
$k = 1 \dots n$ requires \oh{n^2} operations (for each $\boldsymbol{\ell}^{(k)}$,
\oh{k} operations are required to compute $w_k$ and \oh{k} are
required to compute the new entries) and the evaluation of
\eqn{cn} requires \oh{k} operations for each $c_k$, $k=n \dots 0$,
the total cost of the algorithm is in \oh{n^2} operations.

\begin{algorithm}
    \caption{Direct construction of $g_n(x)$}
    \label{alg:direct}
    \begin{algorithmic}[1]
        \STATE $\mathbf p^{(0)} \leftarrow 0$, $\mathbf p^{(1)} \leftarrow \mathbf x$
            \hfill ({\em init $\mathbf P^{(n)}$})
        \FOR{$i=2 \dots n$}
            \STATE $\mathbf p^{(i)} \leftarrow \left( (\mathbf x + \beta_{i-1}) .* \mathbf p^{(i-1)} - \gamma_{i-1} \mathbf p^{(i-2)}\right) / \alpha_{i-1}$
                \hfill ({\em fill $\mathbf P^{(n)}$})
        \ENDFOR
        \STATE $\ell^{(1)}_0 \leftarrow 1$, $w_1 \leftarrow x_1 - x_0$
            \hfill ({\em init $\ell^{(1)}_0$ and $w_1$})\label{alg:direct_ell1}
        \FOR{$i = 2 \dots n$}
            \STATE $w_i \leftarrow 1$ \hfill ({\em construct $w_i$})
            \FOR{$j = 0 \dots i-1$}
                \STATE $w_i \leftarrow w_i (x_i - x_j)$
            \ENDFOR
            \FOR{$j = 0 \dots i-2$}
                \STATE $\ell^{(i)}_j \leftarrow -\frac{\ell^{(i-1)}_j}{x_i - x_j} \frac{w_i}{w_{i-1}}$
                    \hfill ({\em compute the $\ell^{(i)}_j$, \eqn{ellni}})
            \ENDFOR
            \STATE $\ell^{(i)}_{i-1} \leftarrow \frac{w_i}{w_{i-1}(x_i - x_{i-1})}$
                \hfill ({\em compute $\ell^{(i)}_{i-1}$, \eqn{ellnnm1}})
        \ENDFOR \label{alg:direct_ell2}
        \FOR{$i = n \dots 1$} \label{alg:direct_loop1}
            \STATE $c_i \leftarrow \frac{\boldsymbol{\ell}^{(i)} \mathbf f(0:i-1) - f_i}{\boldsymbol{\ell}^{(i)} \mathbf p^{(i-1)}(0:i-1) - p_i(x_i)}$
                \hfill ({\em compute coefficient $c_i$, \eqn{cn}})
            \STATE $\mathbf f \leftarrow \mathbf f - c_i \mathbf p^{(i)}$
                \hfill ({\em update the right-hand side $\mathbf f$})
        \ENDFOR
        \STATE $c_0 \leftarrow f_0 / \mathbf p^{(0)}(0)$
            \hfill ({\em compute the final $c_0$})\label{alg:direct_loop2}
    \end{algorithmic}
\end{algorithm}

Note that, as opposed to the algorithm presented in Section~\ref{sec:first},
this algorithm can not be extended to update or downdate an interpolation.
It has an advantage, however, when multiple right-hand sides, \ie
interpolations over the same set of nodes, are to be computed.
In such a case, the vectors $\boldsymbol{\ell}^{(k)}$, $k=1\dots n$
need to be computed only once (Lines~\ref{alg:direct_ell1} to
\ref{alg:direct_ell2} of Algorithm~\ref{alg:direct}).
For any new vector $\mathbf f$, only the Lines~\ref{alg:direct_loop1}
to \ref{alg:direct_loop2} need to be re-evaluated.

\section{Results}
\label{sec:results}

To assess the stability of the two new interpolation routines described
herein, we will follow the methodology used by \citeasnoun{ref:Higham1988}.
Higham defines a set of interpolations
consisting of all combinations of the nodes
\begin{equation*}
    \begin{array}{lrcll}
        \mbox{A1:} & \quad x_{i} & = & -\cos (i\pi/n), & \mbox{(extrema of $T_n(x)$)} \\
        \mbox{A2:} & x_{i} & = & -\cos \left[ (i + \frac{1}{2})\pi / (n+1) \right], \quad & \mbox{(zeros of $T_{n+1}(x)$)} \\
        \mbox{A3:} & x_i & = & -1 + 2i/n, & \mbox{(equidistant on $[-1,1]$)}\\
        \mbox{A4:} & x_i & = & i/n, & \mbox{(equidistant on $[0,1]$)}
    \end{array}
\end{equation*}
with the right-hand sides
\begin{equation*}
    \begin{array}{lrcl}
        \mbox{F1:} & \quad f_i & = & (-1)^i, \\
        \mbox{F2:} & f & = & ( 1 , 0 , \dots , 0 )^\mathsf{T} \\
        \mbox{F3:} & f_i & = & 1 / (1 + 25x_i^2).
    \end{array}
\end{equation*}
for $i=0 \dots n$.

To avoid instabilities due to unfortunate orderings of the nodes $x_i$,
the nodes and corresponding function values
were re-ordered according to the same permutation that would
be produced by Gaussian elimination with partial pivoting applied
to the Vandermonde-like matrix, as described in \cite{ref:Higham1990}.
This ordering is optimal for the Bj\"orck-Pereyra and Higham algorithms
and produces good results for the two new algorithms described herein.

For each combination of nodes and right-hand sides, we compute, following Higham,
the coefficients $\mathbf c$ for the Chebyshev base (see 
\eqn{coeffs_chebyshev}) for $n=5$, 10, 20 and 30 and compute the quantities
\begin{equation*}
    \mbox{ERR} = \frac{\| \mathbf c - \mathbf c^\star \|_2}{u\|\mathbf c^\star\|_2}, \quad
    \mbox{RES} = \frac{\| \mathbf f - \mathbf P \mathbf c \|_2}{u\|\mathbf c^\star\|_2},
\end{equation*}
where $\mathbf c^\star$ is the exact solution and $u$ is the unit
roundoff\footnote{All results were computed using IEEE 754
double-precision arithmetic and hence $u \approx 2.2\times 10^{-16}$.}
as defined by \citeasnoun[Section 2.4.2]{ref:Golub1996}.
Note that for the special case of A1 or A2 using the Chebyshev
base, the coefficients can also be computed efficiently and reliably
using the Fast Fourier Transform \cite{ref:Battles2004}.

Results were computed using Gaussian elimination ({\small GE}\footnote{For the tests
in this section, Matlab's backslash-operator, which uses partial
pivoting, was used. In cases where the matrix is rank-deficient, a minimum-norm
solution is returned.}), Higham's extension
of the algorithm of Bj\"orck and Pereyra ({\small BP/H}\footnote{Algorithm~1 
in \cite{ref:Higham1988} was implemented in Matlab.}), the incremental
\alg{incr} ({\small INCR}) and
the direct \alg{direct} ({\small DIRECT}).
The exact values were computed in Maple \cite{ref:Maple} with 50 decimal
digits of precision using the {\tt interp} function therein.

Results were also computed for the interpolation downdate ({\small DEL})
described in \alg{downdate}.
Starting from $\mathbf c^\star$ and $\boldsymbol\eta^\star$, the exact
coefficients for the interpolation $g_n(x)$ and the Newton polynomial
$\pi_{n+1}(x)$ respectively, we compute the coefficients $\mathbf c^{(n-1)}$
and $\boldsymbol\eta^{(n)}$ for $g_{n-1}(x)$
and $\pi_{n}(x)$, resulting from the removal of the rightmost function value
$f_k$ at $x_k$, $k=\arg\max_i x_i$.
The exact coefficients $\hat{\mathbf c}^\star$ {\em after} deletion
were computed and used to compute the quantities {\small ERR} and {\small RES}.

\begin{table}
    \begin{small}\begin{center}
    \begin{tabular}{l|cc|cc|cc|cc|cc}
        \multicolumn{1}{c}{} & \multicolumn{2}{c}{{\small GE}} &  \multicolumn{2}{c}{{\small BP/H}} &  \multicolumn{2}{c}{{\small INCR}} &  \multicolumn{2}{c}{{\small DIRECT}} &  \multicolumn{2}{c}{{\small DEL}} \\ 
        $n$ & {\small ERR} & {\small RES} & {\small ERR} & {\small RES} & {\small ERR} & {\small RES} & {\small ERR} & {\small RES} & {\small ERR} & {\small RES} \\ \hline
$5$ & $0.00$ & $0.50$ & $2.00$ & $3.94$ & $\mathbf{3.20}$ & $\mathbf{4.30}$ & $0.00$ & $0.50$ & $0.00$ & $0.72$ \\
$10$ & $3.20$ & $5.32$ & $\mathbf{1.20\mbox{e}1}$ & $\mathbf{2.79\mbox{e}1}$ & $7.76$ & $2.21\mbox{e}1$ & $2.26$ & $7.00$ & $0.00$ & $2.47$ \\
$20$ & $7.28$ & $2.16\mbox{e}1$ & $\mathbf{1.61\mbox{e}2}$ & $\mathbf{5.27\mbox{e}2}$ & $8.92$ & $3.64\mbox{e}1$ & $9.14$ & $2.60\mbox{e}1$ & $0.00$ & $3.40\mbox{e}1$ \\
$30$ & $2.61$ & $1.03\mbox{e}1$ & $\mathbf{6.72\mbox{e}2}$ & $\mathbf{2.65\mbox{e}3}$ & $2.08\mbox{e}1$ & $1.11\mbox{e}2$ & $3.51$ & $1.60\mbox{e}1$ & $0.00$ & $8.76\mbox{e}1$ \\
    \hline \end{tabular}\end{center}
    \caption{Results for problem A1/F1.}
    \label{tab:a1f1}
    \end{small}
\end{table}

\begin{table}
    \begin{small}\begin{center}
    \begin{tabular}{l|cc|cc|cc|cc|cc}
        \multicolumn{1}{c}{} & \multicolumn{2}{c}{{\small GE}} &  \multicolumn{2}{c}{{\small BP/H}} &  \multicolumn{2}{c}{{\small INCR}} &  \multicolumn{2}{c}{{\small DIRECT}} &  \multicolumn{2}{c}{{\small DEL}} \\ 
        $n$ & {\small ERR} & {\small RES} & {\small ERR} & {\small RES} & {\small ERR} & {\small RES} & {\small ERR} & {\small RES} & {\small ERR} & {\small RES} \\ \hline
$5$ & $0.97$ & $2.40$ & $\mathbf{1.82}$ & $\mathbf{3.19}$ & $0.93$ & $2.07$ & $0.88$ & $1.68$ & $0.00$ & $1.70$ \\
$10$ & $2.22$ & $4.19$ & $\mathbf{1.37\mbox{e}1}$ & $\mathbf{3.77\mbox{e}1}$ & $4.94$ & $1.32\mbox{e}1$ & $1.93$ & $3.00$ & $0.00$ & $3.80$ \\
$20$ & $2.11\mbox{e}1$ & $8.96$ & $\mathbf{9.93\mbox{e}1}$ & $\mathbf{3.47\mbox{e}2}$ & $2.24\mbox{e}1$ & $4.80\mbox{e}1$ & $1.80\mbox{e}1$ & $1.27\mbox{e}1$ & $0.00$ & $5.39\mbox{e}1$ \\
$30$ & $3.63\mbox{e}1$ & $1.36\mbox{e}1$ & $\mathbf{1.27\mbox{e}2}$ & $\mathbf{4.84\mbox{e}2}$ & $5.55\mbox{e}1$ & $2.27\mbox{e}2$ & $4.31\mbox{e}1$ & $2.88\mbox{e}1$ & $0.00$ & $1.19\mbox{e}2$ \\
    \hline \end{tabular}\end{center}
    \caption{Results for problem A1/F2.}
    \end{small}
\end{table}

\begin{table}
    \begin{small}\begin{center}
    \begin{tabular}{l|cc|cc|cc|cc|cc}
        \multicolumn{1}{c}{} & \multicolumn{2}{c}{{\small GE}} &  \multicolumn{2}{c}{{\small BP/H}} &  \multicolumn{2}{c}{{\small INCR}} &  \multicolumn{2}{c}{{\small DIRECT}} &  \multicolumn{2}{c}{{\small DEL}} \\ 
        $n$ & {\small ERR} & {\small RES} & {\small ERR} & {\small RES} & {\small ERR} & {\small RES} & {\small ERR} & {\small RES} & {\small ERR} & {\small RES} \\ \hline
$5$ & $1.26$ & $2.67$ & $1.16$ & $2.27$ & $1.26$ & $2.70$ & $\mathbf{1.43}$ & $\mathbf{2.99}$ & $0.00$ & $1.18$ \\
$10$ & $2.12$ & $5.30$ & $\mathbf{7.27}$ & $\mathbf{1.66\mbox{e}1}$ & $1.83$ & $3.63$ & $1.19$ & $1.87$ & $0.32$ & $2.67$ \\
$20$ & $1.13$ & $4.81$ & $\mathbf{8.28}$ & $\mathbf{2.68\mbox{e}1}$ & $1.78$ & $6.74$ & $3.02$ & $1.01\mbox{e}1$ & $0.05$ & $5.25$ \\
$30$ & $1.99$ & $7.35$ & $\mathbf{6.33}$ & $\mathbf{2.53\mbox{e}1}$ & $1.14$ & $5.86$ & $2.55$ & $1.04\mbox{e}1$ & $0.61$ & $6.09$ \\
    \hline \end{tabular}\end{center}
    \caption{Results for problem A1/F3.}
    \end{small}
\end{table}

\begin{table}
    \begin{small}\begin{center}
    \begin{tabular}{l|cc|cc|cc|cc|cc}
        \multicolumn{1}{c}{} & \multicolumn{2}{c}{{\small GE}} &  \multicolumn{2}{c}{{\small BP/H}} &  \multicolumn{2}{c}{{\small INCR}} &  \multicolumn{2}{c}{{\small DIRECT}} &  \multicolumn{2}{c}{{\small DEL}} \\ 
        $n$ & {\small ERR} & {\small RES} & {\small ERR} & {\small RES} & {\small ERR} & {\small RES} & {\small ERR} & {\small RES} & {\small ERR} & {\small RES} \\ \hline
$5$ & $3.55$ & $0.70$ & $\mathbf{5.07}$ & $\mathbf{7.31}$ & $3.79$ & $1.58$ & $4.96$ & $4.96$ & $1.48$ & $3.16$ \\
$10$ & $1.19\mbox{e}1$ & $3.65$ & $\mathbf{2.03\mbox{e}1}$ & $\mathbf{3.31\mbox{e}1}$ & $8.34$ & $1.08\mbox{e}1$ & $1.16\mbox{e}1$ & $2.24$ & $7.23$ & $7.47$ \\
$20$ & $1.66\mbox{e}1$ & $1.23\mbox{e}1$ & $4.64\mbox{e}1$ & $\mathbf{1.63\mbox{e}2}$ & $\mathbf{5.57\mbox{e}1}$ & $1.41\mbox{e}2$ & $1.61\mbox{e}1$ & $2.58$ & $1.69\mbox{e}1$ & $2.41\mbox{e}1$ \\
$30$ & $6.48\mbox{e}1$ & $2.28\mbox{e}1$ & $\mathbf{1.24\mbox{e}2}$ & $\mathbf{4.36\mbox{e}2}$ & $4.45\mbox{e}1$ & $2.29\mbox{e}2$ & $6.49\mbox{e}1$ & $5.98$ & $4.71\mbox{e}1$ & $9.82\mbox{e}1$ \\
    \hline \end{tabular}\end{center}
    \caption{Results for problem A2/F1.}
    \end{small}
\end{table}

\begin{table}
    \begin{small}\begin{center}
    \begin{tabular}{l|cc|cc|cc|cc|cc}
        \multicolumn{1}{c}{} & \multicolumn{2}{c}{{\small GE}} &  \multicolumn{2}{c}{{\small BP/H}} &  \multicolumn{2}{c}{{\small INCR}} &  \multicolumn{2}{c}{{\small DIRECT}} &  \multicolumn{2}{c}{{\small DEL}} \\ 
        $n$ & {\small ERR} & {\small RES} & {\small ERR} & {\small RES} & {\small ERR} & {\small RES} & {\small ERR} & {\small RES} & {\small ERR} & {\small RES} \\ \hline
$5$ & $2.57$ & $3.29$ & $\mathbf{2.88}$ & $\mathbf{4.52}$ & $2.35$ & $0.97$ & $2.30$ & $1.20$ & $0.60$ & $3.25$ \\
$10$ & $5.27$ & $3.89$ & $\mathbf{1.67\mbox{e}1}$ & $\mathbf{4.59\mbox{e}1}$ & $4.94$ & $5.23$ & $3.94$ & $3.75$ & $4.65$ & $1.03\mbox{e}1$ \\
$20$ & $8.40$ & $8.80$ & $\mathbf{5.09\mbox{e}1}$ & $\mathbf{1.63\mbox{e}2}$ & $3.89\mbox{e}1$ & $9.85\mbox{e}1$ & $8.44$ & $4.04$ & $1.15\mbox{e}1$ & $2.43\mbox{e}1$ \\
$30$ & $3.39\mbox{e}1$ & $1.95\mbox{e}1$ & $\mathbf{1.17\mbox{e}2}$ & $\mathbf{4.31\mbox{e}2}$ & $3.00\mbox{e}1$ & $1.99\mbox{e}1$ & $3.42\mbox{e}1$ & $2.55\mbox{e}1$ & $3.27\mbox{e}1$ & $1.13\mbox{e}2$ \\
    \hline \end{tabular}\end{center}
    \caption{Results for problem A2/F2.}
    \end{small}
\end{table}

\begin{table}
    \begin{small}\begin{center}
    \begin{tabular}{l|cc|cc|cc|cc|cc}
        \multicolumn{1}{c}{} & \multicolumn{2}{c}{{\small GE}} &  \multicolumn{2}{c}{{\small BP/H}} &  \multicolumn{2}{c}{{\small INCR}} &  \multicolumn{2}{c}{{\small DIRECT}} &  \multicolumn{2}{c}{{\small DEL}} \\ 
        $n$ & {\small ERR} & {\small RES} & {\small ERR} & {\small RES} & {\small ERR} & {\small RES} & {\small ERR} & {\small RES} & {\small ERR} & {\small RES} \\ \hline
$5$ & $\mathbf{1.44}$ & $0.11$ & $1.40$ & $\mathbf{1.71}$ & $1.17$ & $1.45$ & $1.12$ & $1.67$ & $0.00$ & $2.76$ \\
$10$ & $2.73$ & $3.22$ & $\mathbf{6.06}$ & $\mathbf{1.14\mbox{e}1}$ & $2.86$ & $3.75$ & $3.64$ & $6.34$ & $6.56$ & $7.47$ \\
$20$ & $1.52$ & $5.59$ & $\mathbf{7.34}$ & $\mathbf{2.42\mbox{e}1}$ & $2.06$ & $6.34$ & $3.92$ & $1.34\mbox{e}1$ & $1.27\mbox{e}1$ & $1.81\mbox{e}1$ \\
$30$ & $2.81$ & $1.19\mbox{e}1$ & $\mathbf{7.24}$ & $\mathbf{3.08\mbox{e}1}$ & $1.65$ & $6.07$ & $2.79$ & $9.46$ & $1.41\mbox{e}1$ & $2.94\mbox{e}1$ \\
    \hline \end{tabular}\end{center}
    \caption{Results for problem A2/F3.} \label{tab:a2f3}
    \end{small}
\end{table}

\begin{table}
    \begin{small}\begin{center}
    \begin{tabular}{l|cc|cc|cc|cc|cc}
        \multicolumn{1}{c}{} & \multicolumn{2}{c}{{\small GE}} &  \multicolumn{2}{c}{{\small BP/H}} &  \multicolumn{2}{c}{{\small INCR}} &  \multicolumn{2}{c}{{\small DIRECT}} &  \multicolumn{2}{c}{{\small DEL}} \\ 
        $n$ & {\small ERR} & {\small RES} & {\small ERR} & {\small RES} & {\small ERR} & {\small RES} & {\small ERR} & {\small RES} & {\small ERR} & {\small RES} \\ \hline
$5$ & $1.41$ & $0.58$ & $\mathbf{6.12}$ & $\mathbf{1.01\mbox{e}1}$ & $2.04$ & $0.97$ & $1.48$ & $0.73$ & $0.60$ & $0.78$ \\
$10$ & $2.99$ & $2.50$ & $\mathbf{1.18\mbox{e}1}$ & $\mathbf{2.00\mbox{e}1}$ & $2.16$ & $1.52$ & $2.62$ & $3.11$ & $0.50$ & $1.42$ \\
$20$ & $\mathbf{2.26\mbox{e}3}$ & $4.85$ & $2.01\mbox{e}1$ & $\mathbf{7.19\mbox{e}1}$ & $3.58\mbox{e}1$ & $1.02\mbox{e}1$ & $6.23\mbox{e}1$ & $3.97$ & $0.55$ & $2.77$ \\
$30$ & $\mathbf{1.39\mbox{e}6}$ & $9.51$ & $3.90\mbox{e}1$ & $\mathbf{1.34\mbox{e}2}$ & $5.42\mbox{e}4$ & $2.60\mbox{e}1$ & $3.07\mbox{e}2$ & $7.33$ & $0.55$ & $2.20$ \\
    \hline \end{tabular}\end{center}
    \caption{Results for problem A3/F1.} \label{tab:a3f1}
    \end{small}
\end{table}

\begin{table}
    \begin{small}\begin{center}
    \begin{tabular}{l|cc|cc|cc|cc|cc}
        \multicolumn{1}{c}{} & \multicolumn{2}{c}{{\small GE}} &  \multicolumn{2}{c}{{\small BP/H}} &  \multicolumn{2}{c}{{\small INCR}} &  \multicolumn{2}{c}{{\small DIRECT}} &  \multicolumn{2}{c}{{\small DEL}} \\ 
        $n$ & {\small ERR} & {\small RES} & {\small ERR} & {\small RES} & {\small ERR} & {\small RES} & {\small ERR} & {\small RES} & {\small ERR} & {\small RES} \\ \hline
$5$ & $1.02$ & $1.98$ & $\mathbf{2.19}$ & $\mathbf{3.53}$ & $0.69$ & $0.80$ & $0.73$ & $1.40$ & $0.40$ & $0.86$ \\
$10$ & $\mathbf{4.91}$ & $7.58$ & $3.83$ & $\mathbf{1.03\mbox{e}1}$ & $1.05$ & $1.82$ & $1.61$ & $2.46$ & $0.47$ & $0.85$ \\
$20$ & $\mathbf{3.65\mbox{e}3}$ & $1.07\mbox{e}1$ & $8.32$ & $\mathbf{2.54\mbox{e}1}$ & $5.73\mbox{e}2$ & $2.94$ & $1.29$ & $4.25$ & $0.63$ & $2.26$ \\
$30$ & $\mathbf{4.02\mbox{e}5}$ & $9.93$ & $3.23\mbox{e}1$ & $\mathbf{1.09\mbox{e}2}$ & $1.35\mbox{e}5$ & $1.25\mbox{e}1$ & $4.98$ & $7.32$ & $0.65$ & $3.48$ \\
    \hline \end{tabular}\end{center}
    \caption{Results for problem A3/F2.}
    \end{small}
\end{table}

\begin{table}
    \begin{small}\begin{center}
    \begin{tabular}{l|cc|cc|cc|cc|cc}
        \multicolumn{1}{c}{} & \multicolumn{2}{c}{{\small GE}} &  \multicolumn{2}{c}{{\small BP/H}} &  \multicolumn{2}{c}{{\small INCR}} &  \multicolumn{2}{c}{{\small DIRECT}} &  \multicolumn{2}{c}{{\small DEL}} \\ 
        $n$ & {\small ERR} & {\small RES} & {\small ERR} & {\small RES} & {\small ERR} & {\small RES} & {\small ERR} & {\small RES} & {\small ERR} & {\small RES} \\ \hline
$5$ & $1.48$ & $0.86$ & $\mathbf{1.87}$ & $\mathbf{1.98}$ & $1.78$ & $1.89$ & $1.38$ & $0.36$ & $0.00$ & $2.23$ \\
$10$ & $2.31$ & $3.28$ & $\mathbf{8.96}$ & $\mathbf{2.76\mbox{e}1}$ & $2.00$ & $2.97$ & $3.27$ & $3.80$ & $0.45$ & $1.12$ \\
$20$ & $\mathbf{2.30\mbox{e}3}$ & $5.53$ & $3.81\mbox{e}1$ & $\mathbf{7.38\mbox{e}1}$ & $1.12\mbox{e}2$ & $1.24\mbox{e}1$ & $3.13\mbox{e}1$ & $5.11$ & $0.47$ & $3.39$ \\
$30$ & $\mathbf{1.39\mbox{e}6}$ & $1.10\mbox{e}1$ & $2.28\mbox{e}2$ & $\mathbf{1.88\mbox{e}2}$ & $5.29\mbox{e}4$ & $2.69\mbox{e}1$ & $2.41\mbox{e}2$ & $4.78$ & $0.49$ & $2.29$ \\
    \hline \end{tabular}\end{center}
    \caption{Results for problem A3/F3.} \label{tab:a3f3}
    \end{small}
\end{table}

\begin{table}
    \begin{small}\begin{center}
    \begin{tabular}{l|cc|cc|cc|cc|cc}
        \multicolumn{1}{c}{} & \multicolumn{2}{c}{{\small GE}} &  \multicolumn{2}{c}{{\small BP/H}} &  \multicolumn{2}{c}{{\small INCR}} &  \multicolumn{2}{c}{{\small DIRECT}} &  \multicolumn{2}{c}{{\small DEL}} \\ 
        $n$ & {\small ERR} & {\small RES} & {\small ERR} & {\small RES} & {\small ERR} & {\small RES} & {\small ERR} & {\small RES} & {\small ERR} & {\small RES} \\ \hline
$5$ & $2.01\mbox{e}2$ & $0.55$ & $0.55$ & $\mathbf{0.76}$ & $\mathbf{5.40\mbox{e}2}$ & $0.74$ & $1.57\mbox{e}1$ & $0.48$ & $3.06$ & $3.80$ \\
$10$ & $2.49\mbox{e}4$ & $0.84$ & $0.45$ & $0.76$ & $\mathbf{3.56\mbox{e}6}$ & $\mathbf{1.67}$ & $8.20\mbox{e}2$ & $0.60$ & $1.94$ & $3.55$ \\
$20$ & $\mathbf{2.77\mbox{e}15}$ & $\mathbf{2.53}$ & $1.87$ & $1.69$ & $5.08\mbox{e}14$ & $1.78$ & $6.00\mbox{e}7$ & $1.63$ & $5.93$ & $1.41\mbox{e}1$ \\
$30$ & $4.50\mbox{e}15$ & $0.00$ & $4.44$ & $1.16$ & $     -$ & $     -$ & $3.47\mbox{e}11$ & $1.60$ & $8.93$ & $3.82\mbox{e}1$ \\
    \hline \end{tabular}\end{center}
    \caption{Results for problem A4/F1.}
    \end{small}
\end{table}

\begin{table}
    \begin{small}\begin{center}
    \begin{tabular}{l|cc|cc|cc|cc|cc}
        \multicolumn{1}{c}{} & \multicolumn{2}{c}{{\small GE}} &  \multicolumn{2}{c}{{\small BP/H}} &  \multicolumn{2}{c}{{\small INCR}} &  \multicolumn{2}{c}{{\small DIRECT}} &  \multicolumn{2}{c}{{\small DEL}} \\ 
        $n$ & {\small ERR} & {\small RES} & {\small ERR} & {\small RES} & {\small ERR} & {\small RES} & {\small ERR} & {\small RES} & {\small ERR} & {\small RES} \\ \hline
$5$ & $1.82\mbox{e}2$ & $0.48$ & $0.55$ & $\mathbf{0.65}$ & $\mathbf{3.64\mbox{e}2}$ & $0.37$ & $1.12\mbox{e}1$ & $0.45$ & $1.34$ & $2.27$ \\
$10$ & $1.35\mbox{e}5$ & $0.63$ & $0.40$ & $0.85$ & $\mathbf{3.77\mbox{e}6}$ & $\mathbf{1.55}$ & $6.28\mbox{e}2$ & $0.75$ & $0.99$ & $2.01$ \\
$20$ & $\mathbf{2.81\mbox{e}15}$ & $\mathbf{3.64}$ & $0.71$ & $1.20$ & $6.57\mbox{e}14$ & $0.93$ & $4.05\mbox{e}7$ & $1.57$ & $1.83$ & $3.30$ \\
$30$ & $4.50\mbox{e}15$ & $0.00$ & $0.40$ & $1.27$ & $     -$ & $     -$ & $1.37\mbox{e}11$ & $1.73$ & $2.88$ & $1.27\mbox{e}1$ \\
    \hline \end{tabular}\end{center}
    \caption{Results for problem A4/F2.}
    \end{small}
\end{table}

\begin{table}
    \begin{small}\begin{center}
    \begin{tabular}{l|cc|cc|cc|cc|cc}
        \multicolumn{1}{c}{} & \multicolumn{2}{c}{{\small GE}} &  \multicolumn{2}{c}{{\small BP/H}} &  \multicolumn{2}{c}{{\small INCR}} &  \multicolumn{2}{c}{{\small DIRECT}} &  \multicolumn{2}{c}{{\small DEL}} \\ 
        $n$ & {\small ERR} & {\small RES} & {\small ERR} & {\small RES} & {\small ERR} & {\small RES} & {\small ERR} & {\small RES} & {\small ERR} & {\small RES} \\ \hline
$5$ & $2.41\mbox{e}2$ & $0.23$ & $3.96\mbox{e}1$ & $0.72$ & $\mathbf{8.24\mbox{e}2}$ & $\mathbf{0.73}$ & $8.23$ & $0.62$ & $9.95$ & $1.13\mbox{e}1$ \\
$10$ & $3.79\mbox{e}5$ & $\mathbf{1.38}$ & $1.55\mbox{e}3$ & $0.88$ & $\mathbf{4.94\mbox{e}6}$ & $1.35$ & $3.04\mbox{e}2$ & $0.55$ & $0.26$ & $0.28$ \\
$20$ & $\mathbf{2.81\mbox{e}15}$ & $\mathbf{3.74}$ & $4.84\mbox{e}6$ & $1.21$ & $9.48\mbox{e}14$ & $0.87$ & $3.56\mbox{e}7$ & $1.03$ & $1.45$ & $2.58$ \\
$30$ & $4.50\mbox{e}15$ & $0.00$ & $1.02\mbox{e}11$ & $1.47$ & $     -$ & $     -$ & $1.18\mbox{e}11$ & $1.49$ & $4.56$ & $1.97\mbox{e}1$ \\
    \hline \end{tabular}\end{center}
    \caption{Results for problem A4/F3.}
    \label{tab:a4f3}
    \end{small}
\end{table}

The results are shown in Tables~\ref{tab:a1f1} to \ref{tab:a4f3}. 
For each $n$, the largest values for {\small ERR} and {\small RES}
are highlighted.
For the problem sets over the nodes A1 and A2 (Tables~\ref{tab:a1f1} to \ref{tab:a2f3}),
the condition of the Vandermonde-like
matrix is always $\leq 2$ \cite{ref:Gautschi1983}, resulting in very small errors
for Gaussian elimination.
The Bj\"orck-Pereyra/Higham algorithm generates slightly larger residuals than both the
incremental and direct algorithms for both sets of nodes.
The values for {\small ERR}, however, are usually within the same
order of magnitude for the three algorithms.

For the nodes A3 (Tables~\ref{tab:a3f1} to \ref{tab:a3f3}),
the condition number of the Vandermonde-like matrix is $5.11\mbox{e}6$ for
$n=30$, resulting in the errors of approximately that magnitude
when Gaussian elimination is used.
In general, both the Bj\"orck-Pereyra/Higham and the direct algorithm
generate smaller errors and residues than Gaussian elimination.
The errors for the incremental algorithm are due to
cancellation while evaluating $g_n(x_{n+1})$ for \eqn{anp1} since
the intermediate coefficients $\mathbf c^{(k)}$ are several orders
of magnitude larger than the result\footnote{
In \cite{ref:Higham1988}, Higham shows that the coefficients
can be written as the weighted sum of any of the intermediate coefficients
$c^{(n)}_i = \sum_j \mu_j c^{(k)}_j$, where the $\mu_j$ depend only
on the nodes and the coefficients of the three-term recurrence relation.
If the $\mu_j$ are \oh{1} and the intermediate $c^{(k)}_j$ are much 
larger than the $c^{(n)}_i$, then cancellation is likely to occur in 
the above sum.}.

Finally, the condition number of the Vandermonde-like matrix for the
nodes A4 is $4.26\mbox{e}16$ for $n=30$, making it numerically singular and thus
resulting in the complete failure of Gaussian elimination.
Note that since in such cases Matlab's backslash-operator computes the
minimum norm solution, the resulting residual error {\small RES} is quite small.
For the first two right-hand sides F1 and F2, the Bj\"orck-Pereyra/Higham
algorithm performs significantly better than the two new algorithms,
since the magnitude of the intermediate coefficients does not vary significantly.
For the right-hand side F3, however, the errors are larger,
caused by truncation in computing the Newton coefficients $a_i$.
The incremental algorithm fails completely for all right-hand sides
since the intermediate and final coefficients $\mathbf c^{(k)}$, $k \leq n$, are
more than ten orders of magnitude larger than the function values\footnote{
$\|\mathbf c^\star\| = 2.23\mbox{e}13$ for F3 and $n=30$.} and
the numerical condition of $g_n(x_{n+1})$ in \eqn{anp1} exceeds
machine precision, resulting in numerical overflow.
These relatively large coefficients also cause problems for the direct
algorithm when evaluating the right-hand side of \eqn{second_new}, where
the original function values are clobbered by the subtraction of the
much larger $\mathbf p^{(n)}c_n$.

The errors and residuals for the downdate algorithm are shown in the 
rightmost columns of Tables~\ref{tab:a1f1} to \ref{tab:a4f3}.
In general the errors of the downdate are relatively small
for all test cases.
The larger residues, \eg for A2/F2, are due to cancellation
in the final subtraction in \alg{downdate}, Line~\ref{alg:downdate_last}.

\section{Conclusions}

We have presented here two new algorithms for the construction of
polynomial interpolations.
The first algorithm (Algorithm~\ref{alg:incr}) offers no substantial
improvement over that of Bj\"orck-Pereyra/Higham except that it can
be easily downdated.
The second algorithm, which does not allow for updates nor downdates,
is slightly more stable than the other algorithms tested and is more
efficient when multiple right-hand sides need to be computed over the
same set of nodes.

\bibliographystyle{kluwer}
\bibliography{interp}

\end{document}